\documentclass[twoside,11pt]{article}
\usepackage[dvipsnames]{xcolor}

\usepackage{custom}

\usepackage{algorithm, algorithmicx}
\usepackage{algpseudocode}
\usepackage{subfig}
\usepackage{physics}
\usepackage{subcaption}
\usepackage{float}
\usepackage{graphicx}
\usepackage{adjustbox}

\usepackage{amsfonts}
\usepackage{amsmath}
\usepackage{amssymb}
\usepackage{mathtools}

\usepackage{pstricks, pst-plot}	
\usepackage[figuresright]{rotating}

\usepackage[english]{babel}

\usepackage{blindtext}
\usepackage[hypcap=true,labelfont={bf},labelsep=period]{caption}
\usepackage{braket}

\usepackage{booktabs}

\usepackage{gensymb}

\title{Global Optimization via Quadratic Disjunctive Programming for Water Networks Design with Energy Recovery}

\author[1]{Carolina Trist\'an}
\author[2]{Marcos Fallanza}
\author[2]{Raquel Ib\'añez}
\author[3]{Ignacio E. Grossmann}
\author[1,4,5]{David~E.~Bernal~Neira }

\affil[1]{\RaggedRight Davidson School of Chemical Engineering, Purdue University, West Lafayette, IN, USA}
\affil[2]{Department of Chemical and Biomolecular Engineering, University of Cantabria, Santander, Spain}
\affil[3]{Department of Chemical Engineering, Carnegie Mellon University, Pittsburgh, PA, USA}
\affil[4]{USRA Research Institute for Advanced Computer Science, Mountain View, CA, USA}
\affil[5]{NASA Ames Research Center, Moffett Field, CA, USA}

\email{dbernaln@purdue.edu}

\begin{document}
\justifying
\setlength{\parindent}{0pt}
\maketitle
\thispagestyle{empty}  

\begin{abstract}
Generalized disjunctive programming (GDP) models with bilinear and concave constraints, often seen in water network design, are challenging optimization problems. This work proposes quadratic and piecewise linear approximations for nonlinear terms to reformulate GDP models into quadratic GDP (QGDP) models that suitable solvers may solve more efficiently. We illustrate the benefits of the quadratic reformulation with a water treatment network design problem in which nonconvexities arise from bilinear terms in the mixers’ mass balances and concave investment cost functions of treatment units. Given the similarities with water network design problems, we suggest quadratic approximation for the GDP model for the optimal design of a large-scale reverse electrodialysis (RED) process. This power technology can recover energy from salinity differences between by-product streams of the water sector, such as desalination brine mixed with regenerated wastewater effluents. The solver Gurobi excels in handling QGDP problems, but weighing the problem’s precision and tractability balance is crucial. The piecewise linear approximation yields more accurate, yet larger QGDP models that may require longer optimization times in large-scale process synthesis problems.
\end{abstract}


\section{Introduction} 
\label{sec:intro}
Water network design often involves bilinear mass balances and nonlinear investment cost functions that lead to nonconvex generalized disjunctive programming (GDP) with bilinear and concave constraints \cite{Ruiz2016}. This work proposes replacing the nonlinear equations with piecewise or quadratic approximations to define a quadratic generalized disjunctive program (QGDP) for GDP problems. These problems feature bilinearities that commercial solvers, such as Gurobi, can efficiently solve to global optimality \cite{Achkar2023}.

The general form of a process superstructure GDP optimization model reads as follows:
\begin{equation}\label{eq:gdp}
    \begin{aligned} 
        \min_{x,Y}\ z = &\ f(x) \\
        \text{s.t.} \quad &\ g(x) \leq 0\\
        & \left[
            \begin{gathered}
            Y_{u} \\
            h_{u}(x)\leq 0\\
            \end{gathered}
            \right]
            \vee \left[
            \begin{gathered}
            \neg \ Y_{u} \\
            B^{u}\ x = 0\\
            \end{gathered}
        \right] \quad \forall \ u \in U\\
        &\ \Omega(Y_1,Y_2,\dots,Y_{\mid U \mid}) = True \\
        &\ x \in X \subseteq \mathbb{R}^n\\
        &\ Y_u \in \{True, False\}
    \end{aligned}   
\end{equation}
where continuous bounded variables $x$ optimize an objective function $f(x)$ subject to global constraints, $g(x) \leq 0$, and a set of disjunctions that determines whether unit $u$ is selected or not. Boolean variables, $Y_u$, in each disjunct activate constraints $h_u(x) \leq 0$ relative to that unit when it is selected ($Y_u=True$); otherwise, ($\neg Y_u$) ignore equations in the inactive disjunct and fix some variables to zero $B^u  x = 0$, involved in the inactive unit. The logic constraints ($\Omega(Y)=True$) set conditions for selecting specific units.

Water networks involve water-using process and treatment units, offering several water integration alternatives that reduce freshwater consumption and wastewater generation while minimizing the total network cost subject to a specified discharge limit \cite{Ahmetovic2017}. Reverse electrodialysis (RED) can recover the salinity gradient energy (SGE) embedded in the mixing of the water network’s streams of different salinities, providing a sustainable supply of clean, base-load electricity to regeneration or water-supply processes \cite{Rani2022}. In previous work \cite{Tristan2023}, we developed a GDP model of the RED process that incorporates a detailed model of the RED unit. This GDP resembles the water network design model. Bilinear mass balances in the mixers and nonlinearities in the RED unit model give rise to a nonlinear GDP model. To solve to global optimality, we apply the Global Logic-based Outer Approximation (GLOA) algorithm \cite{Chen2022}, which decomposes the solution to the GDP into a sequence of mixed-integer linear programming (MILP) master problems and reduced nonlinear-programming (NLP) subproblems that need to be solved global optimality to guarantee a global optimum convergence. We explore how the GDP reformulation into a QGDP approach affects the tractability and fidelity of the optimization model, assessing the computational time and accuracy of the optimal solutions for two case studies: water treatment network design and a large-scale RED process synthesis problem.
\section{Case Study 1. Water Treatment Network Design}
\label{sec:case1}
\subsection{Problem Statement}
Given is a set of water streams with known concentrations of contaminants and flow rate. The objective is to find the set of treatment units and interconnections that minimize the cost of the WTN while satisfying maximum concentrations of contaminants in the reclaimed outlet stream. The WTN superstructure consists of a set of treatment units $t\in TU$, contaminated feed streams $f\in FU$ carrying a set of contaminants $j\in J$, and a discharge unit $disch$. The fouled feed waters can be allocated to one or more treatment units or disposed of in the sink unit. Upon treatment, the reclaimed streams can be recycled, forwarded to other treatment units, or discharged into the sink unit.
\subsection{Model formulation}
The WTN design GDP problem reads as \eqref{eq:gdp}. The mass balances are defined in terms of total flows, $F_s$, and contaminants concentration, $C_{j,s}$  \cite{Karuppiah2006,Quesada1995}. Together with the treatment unit cost $CTU_t$, these determine the continuous variables of the problem, i.e., $x=\left[F_s,C_{j,s},CTU_t\right]$. The treatment units have fixed recoveries, $\alpha_{j,t}$ is the recovery of contaminant $j$ in treatment unit $t$; $C_{j,mt,t}$ and $C_{j,t,st}$ are contaminant concentration at the inlet and the outlet of $t$. The treatment unit cost $CTU_t$ consists of a linear operating cost term and a concave capital cost term. Nonconvexities arise from bilinear terms “flows times concentration” in the mixers’ mass balances and concave investment cost functions of treatment units.

Following Eq. \eqref{eq:gdp}, the objective is to minimize the total treatment unit costs $f(x) = \sum_{t\in TU}CTU_t$, the equations set for inexisting units become $B^t x = \left[\sum_{i\in S_{mt}}F_i ,CTU_t \right] = 0\quad \forall \ t \in TU$, and the model constraints become
\begin{align}\label{eq:gh_wtn}
    &g(x) \leq 0 \quad \rightarrow \quad
    \left\{
    \begin{aligned}
        &F_k \ C_{j,k} = \sum_{i\in S_i\subseteq S}F_i \ C_{j,i} \quad \forall \ j \in J, \ k \in S_k\subseteq S \\
        &F_k = \sum_{i\in S_i\subseteq S}F_i \quad \forall \ k \in S_k\subseteq S \\
        &C_{j,k} = C_{j,i} \quad \forall \ j \in J,\ k \in S_k\subseteq S, \ i \in S_i\subseteq S \\
        &F_{i} = \sum_{k \in S_k\subseteq S} F_k \quad \forall \ i \in S_i\subseteq S \\
        &\sum_{i\in S_{dm}\subseteq S}F_i \ C_{j,i} \leq T_{j} \quad \forall \ j\in J
    \end{aligned}
    \right. , \\
    &h_{t}(x) \leq 0 \quad \rightarrow \quad
    \left\{
    \begin{aligned}
        &C_{j,t,st}=(1-\alpha_{j,t})\ C_{j,mt,t} \quad \forall \ j \in J \\
        &F_{mt,t} \ C_{j,mt,t} = \sum_{i\in S_i\subseteq S}F_i \ C_{j,i} \quad \forall \ j \in J \\
        &F_{mt,t} = \sum_{i\in S_i\subseteq S}F_i \\
        &C_{j,t,st} = C_{j,i} \quad \forall \ j \in J,\ i \in S_i\subseteq S \\
        &F_{t,st} = \sum_{k \in S_{st}\subseteq S} F_k \\
        &F_{mt,t} = F_{t,st} \\
        &F_{mt,t} \geq L_t \\
        &CTU_t=\beta_t \ F_{mt,t}+\gamma_t+\theta_t \ F^{0.7}_{mt,t}
    \end{aligned}
    \quad \forall \ t \in TU
    \right.
\end{align}
\section{Case Study 2. RED Process Design}
\label{sec:case2}
\subsection{Problem Statement}
Given is a set of candidate RED units and the high- and low-salinity feeds’ concentration, flow rate, and temperature. The objective is to determine the hydraulic topology and operational conditions of the RED units that maximize the net present value (NPV) of the process. The superstructure definition and notation are outlined in \cite{Tristan2023}.
\subsection{Model Formulation}
From Eq. \eqref{eq:gdp}, $f(x)$ maximizes the NPV of the RED process. The variables $x$ are the molar concentration and flow rate of the streams and the internal variables of the active RED units. The decision variables are the active RED stacks’ operating conditions. The global constraints, $g(x)\leq 0$, specify mass balances that must hold for any selection of alternatives. The set of disjunctions determines whether the RED unit $r$ is active, enforcing $h_r (x)\leq 0$, the RED unit model equations. These factor in the capital and operating costs and set bounds on the RED unit’s internal variables and the inlet and outlet streams’ concentration and flow rate. Nonconvexities arise from bilinear mass balances in the mixers and the RED unit model, the concave investment cost of pumps, and the Nernst electric potential and gross power equations in the RED unit model.
\section{Results and Discussion}
\label{sec:results}
We implement the GDP models using the algebraic modeling language Pyomo \cite{Hart2017} and the Pyomo.GDP extension \cite{Chen2022} for logic-based modeling and optimization. We solve the problems on a Windows 10 (x64) machine with 6 cores processor (Intel\textsuperscript{®} Core™ i7-8700 CPU @3.2 GHz) and 16 GB of RAM using the solver versions from GAMS 34.1.0, and setting a 0.01 \% optimality gap and 3600 s time limit. We use two proposed reformulation strategies for the nonlinear terms in Eq. (1): (a) Problem (QGDP-q): fitting a quadratic function, i.e., $f(x)\approx x' \ Q\ x$ with $Q$ potentially nonconvex, using the \texttt{curve\_fit()} function from the \texttt{SciPy}. (b) Problem (QGDP-pwl): using a piecewise linear approximation using the \texttt{Piecewiselinear} Pyomo functionality. We opted for the incremental model for the piecewise function \cite{Vielma2010} and split the domain into 101 segments which adds binary variables to the GDP problem.
\subsection{Water Treatment Network Design}
The WTN comprises five inlet streams with four contaminants and four treatment units. The contaminant concentration and flow rate of the feed streams, contaminant recovery rates ($\alpha_{j,t}$), minimum flow rate ($L_t$) and cost coefficients ($\beta_t,\ \gamma_t,\ \theta_t$) of the treatment units, and the upper limit on the molar flow of contaminant j in the purified stream ($T_j$), are reported in \cite{MINLP24}. To solve the WTN problem, we convert the GDP models into MINLPs using the Big-M reformulation \cite{Chen2022}. The MINLP model with the concave investment cost is solved with BARON and the MINLPs with the quadratic and piecewise linear approximations with Gurobi and BARON.
\begin{table}[htbp]
  \centering
  \caption{Model size of the WTN GDP with the concave cost term (original) and the quadratic and piecewise linear (pwl) QGDPs, and MINLP solvers’ computational time, minimum WTN cost, and relative error between QGDP and GDP optimal cost}
    \begin{tabular}{lccccc}
    \toprule
    \multicolumn{1}{c}{} & \textbf{Original GDP} & \multicolumn{2}{c}{\textbf{Quadratic QGDP-q}} & \multicolumn{2}{c}{\textbf{pwl QGDP-pwl}} \\
    \midrule
    \# cont. vars  & \multicolumn{1}{c}{239} & \multicolumn{2}{c}{239} & \multicolumn{2}{c}{749} \\
    \# binary vars & \multicolumn{1}{c}{10} & \multicolumn{2}{c}{10} & \multicolumn{2}{c}{510} \\
    \# const (nl) & 329 (33) & \multicolumn{2}{c}{329 (33)} & \multicolumn{2}{c}{1339 (28)} \\
    Solver & BARON & BARON & Gurobi & BARON & Gurobi \\
    CPU time [s] & \multicolumn{1}{c}{16.12} & 16.86 & 6.44  & 2543  & 10.06 \\
    Objective & \$348,337 & \$349,556 & \$349,562 & \$348,337 & \$348,337 \\
    Relative error [\%] & n.a.  & 0.35 & 0.35 & 0.00 & 0.00 \\
    \bottomrule
    \end{tabular}%
  \label{tab:wtn}%
\end{table}%

The QGDP models and the GDP with the concave capital cost term yield the same optimal WTN design, in which all the polluted streams but one are treated in units one and four \cite{MINLP24}. Table \ref{tab:wtn} compares the size, computational time, and optimal solution of the WTN design problem for each modeling approach and solution strategy. With a relative error of 0.35 \%, Gurobi yields an optimal QGDP-q solution in just half the time it takes BARON to solve the original GDP model with the same number of variables and constraints. BARON requires more time than Gurobi to find the optimal solution of the QGDP-q model. Under the pwl approximation’s better fit to the concave capital cost term, Gurobi and BARON obtain the same global optimum to the QGDP-pwl model as BARON does for the GDP. While the pwl approximation offers more accurate models, the increase in size may also render them intractable. For instance, BARON requires almost an hour to find the optimal WTN design. By contrast, Gurobi solves the instance in about the same time as for the QGDP-q.
\subsection{RED Process Design}
The process consists of a set of candidate RED units, which draw energy from the effluent of a real desalination plant that rejects 733 m\textsuperscript{3}/h of brine (1.67 M NaCl, 20 \degree C) \cite{Tristan2020}. A wastewater treatment plant provides an equal volume of low-salinity feedwater (20 mM NaCl) \cite{PerezTalavera2001} for SGE conversion.

Prior work \cite{Tristan2023} showed that the most time-consuming steps in the GLOA algorithm involve the discretized RED unit model, i.e., the initial linearization of the GDP and resolution of the NLP subproblems. As the number of candidate RED units in the process grows, which is expected in large-scale systems, these steps become even more expensive. Therefore, approximating it as a QP model could accelerate the GDP solution without sacrificing fidelity. We compare the solution that maximize the net power from the QP-q and QP-pwl models with the NLP model (with Nernst potential and gross power nonlinear equations) to appraise the QP models’ fidelity and tractability in Table \ref{tab:red}.
\begin{table}[htbp]
  \centering
  \caption{Model size of the RED unit’s NLP without reformulations (original), quadratic and piecewise linear (pwl) QPs, and solver’s computational time, maximum net power, and relative error between QPs and NLP optimal solution.}
    \begin{tabular}{lccccc}
    \toprule
    \multicolumn{1}{c}{} & \textbf{Original NLP} & \multicolumn{2}{c}{\textbf{Quadratic QNLP-q}} & \multicolumn{2}{c}{\textbf{pwl QNLP-pwl}} \\
    \midrule
    \# cont. vars  & \multicolumn{1}{c}{179} & \multicolumn{2}{c}{180} & \multicolumn{2}{c}{1608} \\
    \# binary vars & \multicolumn{1}{c}{n.a.} & \multicolumn{2}{c}{n.a.} & \multicolumn{2}{c}{1400} \\
    \# const (nl) & 182 (101) & \multicolumn{2}{c}{183 (103)} & \multicolumn{2}{c}{3011 (89)} \\
    Solver & BARON & BARON & Gurobi & BARON & Gurobi \\
    CPU time [s] & 3600 & 306 & 115 & 3600 & 1986 \\
    Objective [kW] & 1.013 & 1.043 & 1.043 & 1.008 & 1.016 \\
    Relative error [\%] & n.a.  & 2.9615 & 2.9615 & 0.4936 & 0.2961 \\
    \bottomrule
    \end{tabular}%
  \label{tab:red}%
\end{table}%

The quadratic approximation’s inaccurate fitting of natural logarithms within the Nernst potential leads to errors in the optimal solution for the NLP (Table \ref{tab:red}). Despite this, the computational time is decreased with Gurobi and BARON, with Gurobi outperforming, while keeping almost the same number of variables and constraints. The piecewise linear approximation provides a better fit but significantly increases the size. After an hour, BARON could only reach a suboptimal solution to the QP-pwl with an 27~\% optimality gap. Despite the increase in the size of the QP-pwl, Gurobi even outperforms BARON solving the rigorous NLP model (original NLP), which finds the optimal net power with a 24~\% optimality gap within an hour. The outcomes of this study provide the drive to transform the RED process GDP model into a quadratic one, a subject of future work.
\section{Conclusions}
\label{sec:conclusions}
Water network design problems generally involve bilinearities and concave functions that may lead to multiple local optima. Solving global optimality thereby requires computationally expensive global optimization approaches. This work presents quadratic and piecewise linear approximations for nonlinearities in GDP models involving bilinearities to derive a quadratic GDP model that solvers like Gurobi can efficiently solve. The implications are far-reaching for problems where bilinear constraints are present with constraints that fit well with quadratic functions or piecewise linear (or even piecewise quadratic) approximations. For handling this type of problem, Gurobi is a powerful solver, though carefully considering the trade-off between fidelity and tractability is paramount. Using piecewise linear approximations grants higher accuracy, but also introduces additional variables and constraints to the original problem, making it challenging to find the global optimum within an acceptable timeframe, especially for larger-scale problems that may require advanced decomposition strategies.
\section*{Acknowledgements}
\begin{sloppypar}
Projects TED2021-129874B-I00, PDC2021-120786-I00, and PID2020-115409RB-I00 through EU NextGenerationEU/PRTR and MCIN/AEI/ 10.13039/501100011033. C.T. acknowledges \mbox{PRE2018-086454} fellowship, "ESF Investing in your future", and, with D.B.N., the School of Chemical Engineering at Purdue University’s startup grant. D.B.N. acknowledges the NASA Academic Mission Services, Contract No.~NNA16BD14C.
\end{sloppypar}
\bibliography{References}
\end{document}